\documentclass[12pt, leqno]{amsart}
\usepackage{amsmath}
\usepackage{amssymb}

\usepackage[foot]{amsaddr}

\usepackage{algorithm}
\usepackage[noend]{algpseudocode}

\usepackage{bm}
\usepackage{graphicx}
\usepackage{caption}
\usepackage{subcaption}
\usepackage[hidelinks]{hyperref}

\usepackage{enumerate} 
\usepackage{parskip}
\usepackage{xcolor}
\allowdisplaybreaks

\def\red{\textcolor{red}}

\newtheorem{theorem}{Theorem}[section]

\theoremstyle{definition}

\theoremstyle{remark}
\newtheorem{remark}[theorem]{Remark}
\numberwithin{equation}{section}

\newcommand{\R}{\mathbb{R}}

\DeclareMathOperator*{\argmin}{arg\!\min}

\numberwithin{equation}{section}
\renewcommand{\theequation}{\arabic{section}.\arabic{equation}}

\begin{document}
\title[AA gradient methods with energy for optim
probs] 
{Anderson acceleration of gradient methods with energy for optimization problems}


\author[H. Liu]{Hailiang Liu\textsuperscript{\textdagger}\textsuperscript{\textasteriskcentered}}
\thanks{\textsuperscript{\textdagger}Department of Mathematics, Iowa State University, Ames, IA, USA}

\author[J.H. He]{Jia-Hao He\textsuperscript{\textdaggerdbl}}
\thanks{\textsuperscript{\textdaggerdbl}Department of Agricultural and Biosystems Engineering, Iowa State University, Ames, IA, USA}

\author[X. Tian]{Xuping Tian\textsuperscript{\textdagger}}

\thanks{\textsuperscript{\textasteriskcentered}corresponding author: \url{hliu@iastate.edu}}

\dedicatory{Dedicated to Professor Stanley Osher on the occasion of his 80th birthday}




\begin{abstract}
Anderson acceleration (AA) as an efficient technique for speeding up the convergence of fixed-point iterations may be designed for accelerating an optimization method. 
We propose a novel optimization algorithm by adapting Anderson acceleration to the energy adaptive gradient method (AEGD) [\href{https://arxiv.org/abs/2010.05109}{arXiv:2010.05109}]. The feasibility of our algorithm is examined in light of convergence results for AEGD, though it is not a fixed-point iteration. 
We also quantify the accelerated convergence rate of AA for gradient descent by a factor of the gain at each implementation of the Anderson mixing. Our experimental results show that the proposed algorithm requires little tuning of hyperparameters and exhibits superior fast convergence.
\end{abstract}

\maketitle


\section{Introduction}

We are concerned with nonlinear acceleration of gradient methods with energy for optimization problems of form  
\begin{equation}\label{probsetup}
    \min_{x\in \mathbb{R}^n} f(x),  
\end{equation}
where $f:\mathbb{R}^n \rightarrow \mathbb{R}$ is assumed to be differentiable and bounded from below. Gradient descent (GD), known for its easy implementation and low computational cost, has been one of the most popular optimization algorithms. However, the practical performance of this well-known first-order method is often limited to convex problems. For non-convex or ill-conditioned problems, GD can develop a zig-zag pattern of subsequent iterates as iterations progress, resulting in slow convergence.
Multiple acceleration techniques have been proposed to address these deficiencies,  
including those in the fast gradient methods \cite{Nest1, Nest2} and the momentum or heavy-ball method \cite{momentum}. Although reasonably effective and computationally efficient, the gradient descent method with these techniques might
still suffer from the step size limitation. 

AEGD (adaptive gradient descent with energy) is an improved variant of GD originally introduced in  \cite{LT20} to address the step size issue. 
AEGD adjusts the effective step size by a transformed gradient and an energy variable. The method, when applied to the problem \eqref{probsetup},   includes two ingredients: the base update rule:
\begin{equation}\label{aegd0}
x_{k+1}=x_k +2\eta r_{k+1}v_k, \quad r_{k+1} = \frac{r_{k}}{1+2\eta v^2_{k}},
\end{equation} 
and the evaluation of the transformed gradient $v_k$ as 
\begin{equation}\label{vt}
v_{k}=\frac{\nabla f(x_k)}{2\sqrt{f(x_k)+c}}.
\end{equation}
AEGD is unconditionally energy stable with guaranteed convergence in energy regardless of the size of the base learning rate $\eta>0$ and how $v_k$ is evaluated. Adding momentum through modifying $v_k$ can further accelerate the convergence in both deterministic and stochastic settings \cite{AEGDM, ODE, SGEM}. 

In this paper, we aim to speed up AEGD with another acceleration technique, called Anderson acceleration (AA) \cite{Anderson}.
As a nonlinear acceleration technique, AA has been used for decades to speed up the convergence of a sequence by using previous iterates to extrapolate a new estimate. AA is widely known in the context  of fixed-point iterations \cite{AAFP, ZDB20}. Recently, AA has also been adapted to various optimization algorithms 
\cite{AA-PGD, BM21, AA-optim1, AA-optim2, Regular2, AA-optim3}; More specifically, we refer to \cite{AA-PGD} for the proximal gradient algorithm, \cite{BM21} for the coordinate descent, \cite{AA-optim1} for primal-dual methods, \cite{AA-optim2} for the ADMM method, \cite{Regular2} for the Douglas-Rachford algorithm, and \cite{AA-optim3} for stochastic algorithms. All these algorithms can be formulated as involving fixed-point iterations. In contrast, our work applies AA to an optimization algorithm that is not a fixed-point iteration.

Convergence analysis of AA for general fixed-point maps has been reported in the literature only recently \cite{ConvAnal, EDIIS, TE17, BC21}. 
If the fixed point map is a contraction, the first local convergence results for AA are provided in \cite{ConvAnal, EDIIS} for either $m=1$ or general $m$ if the coefficients in the linear combination remain bounded; see also \cite{TE17, BC21} for further extensions to nonsmooth or inaccurate function evaluations.  
But these results do not show how much AA can help accelerate the convergence rate. Our convergence rates on AA with GD are for general $m>1$,  showing both the linear rate of the GD for strongly convex objectives and the additional gain from AA. The amount of acceleration at each AA implementation is quantified by using a projection operator. This improves a convergence bound in \cite{AA-optim4} (see Remark \ref{rem3.3}).  For contracting fixed point maps, we refer to \cite{EP20} for a refined estimate of the convergence gain by AA, and \cite{AA-nest, WS21} for further asymptotic linear convergence speed of AA when applied to more general fixed-point iteration and ADMM, respectively.  
Sharper local convergence results 
of AA remain a hot research topic in this area.

The motivation of this work is to incorporate the two acceleration techniques: adaptive gradient descent with energy and AA, into one framework so that the resulting method can take advantage of both techniques. We summarize the main contributions of this work as follows:
\begin{itemize}
\item With AA for GD, we quantify the gain of convergence rate for both quadratic and non-quadratic minimization problems by a ratio factor of the optimized gradient compared to that of the vanilla GD.
\item We present a modification of AA for AEGD by introducing an extra hyperparameter that controls the implementation frequency of AA. Our algorithm (Algorithm \ref{aaegd}) is shown to be insensitive to both the step size $\eta$ and the length of the AA window $m$, thus requiring little hyperparameter tuning.
\item We also adapt our algorithm for the proximal gradient algorithm (Algorithm \ref{aapeg}) for solving a class of composite optimization problems whose objective function is given by the summation of a general smooth and nonsmooth component.
\item We verify the superior performance of our  algorithms on both convex and non-convex tests, also on two machine learning problems. In all cases, our proposed algorithms exhibit faster convergence than AA with vanilla gradient descent.
\end{itemize}

\subsection{Related work} As its name suggests, AA or Anderson mixing (AM) is an acceleration algorithm originally developed by D. G. Anderson \cite{Anderson} for solving some  nonlinear integral equations, and it is now commonly used in scientific computations for problems that can be regarded as involving a fixed-point iteration. Applications include flow problems \cite{flow1, flow2, GRVZ22}, electronic structure computations \cite{Anderson, FS09}, 
wave propagation \cite{Ya21, LL22}, and deep learning \cite{Pasi21, RL18}. In a broad literature on this subject, these papers are only meant to show the variety of results obtained by the AA. 
The connection between AA and some other classical methods also facilitates our understanding of AA. As originally noted by Eyert \cite{quasinewton} and clarified by  Fang and Saad \cite{FS09}, AA is remarkably related to a multisecant quasi-Newton updating. For linear iterations, AA with full-memory (i.e., $m=\infty$ in Algorithm \ref{aagd}) is essentially equivalent to the generalized minimal residual (GMRES), as shown in \cite{AAFP, OW00, PE13}.
The acceleration property of AA is also understood from its close relation to Pulay mixing \cite{Pu80} and to DIIS (Direct Inversion on the Iterative Subsapce) \cite{DIIS, EDIIS}.  AA is also related to other acceleration techniques for vector sequences, where they are known under various names, e.g. minimal polynomial extrapolation \cite{SF87} or  reduced rank extrapolation \cite{E79}.

\subsection{Organization}
The paper is organized as follows. In Section 2, we start by reviewing the sequence acceleration method -- AA. In Section 3, we recollect the theoretical results of the fixed-point iterations with AA and then extend the result to GD with nonlinear iterations. In Section 4, we recall the convergence of AEGD to explain how to implement AA with AEGD, and further introduce AA-AEGD to the application of proximal gradient descent. Finally, we present numerical results in Section 5 and the conclusion in Section 6.

\section{Review of Anderson acceleration} 

AA is an efficient acceleration method for fixed-point iterations $x_{k+1}=G(x_{k})$, which is assumed to converge to $x^*$ such that $x^*=G(x^*)$. The key idea of the AA  is to form a new extrapolation point by making better use of past iterates. To illustrate this method, we consider a window of sequence of $k+1$  elements and their updates, denoted by $X=[x_0, x_1, \cdots x_{k}]$ and $Y=G(X)=[G(x_0), \cdots, G(x_{k})]$, respectively. And denote $R_k=G(x_k)-x_k$ as the residual at $x_{k}$. The goal is to find 
$\alpha \in \mathbb{R}^{k+1}$ such that 
\begin{equation}\label{aay}
y_{k+1}=X\alpha ,  \quad \alpha^\top \mathbf{1}=1    
\end{equation}
has a faster convergence rate in the sense that 
$$
\|y_{k+1}-x^*\|<< \|x_{k+1}-x^*\|, \quad x_{k+1}=G(x_{k}).
$$
AA determines $\alpha$ by solving the following minimization problem
\begin{equation}\label{a_prev}
    \alpha={\rm argmin}_{\alpha^\top \mathbf{1}=1}\|R\alpha\|_2,
\end{equation}
where $R=[R_0,\cdots, R_k]$ is the residual matrix. As an extension to AA, one may modify (\ref{aay}) by introducing a variable relaxation parameter $\beta \in (0, 1]$:
$$
y_{k+1}=[(1-\beta)X +\beta Y]\alpha.
$$
We summarize the method with the window of length $m$ in Algorithm \ref{aagd}. 
\begin{algorithm}
\caption{Anderson Acceleration (AA)} 
\label{aagd}
\begin{algorithmic}[1] 
\Require $G(\cdot)$: fixed iteration map; $x_0$: initialization of the parameter; $\eta$: step size, $m$: length of the window; $\beta$: relaxation parameter; $K$: total number of iterations.
\State Set $x_1 = G(x_0)$ and $R_0=G(x_0)-x_0$.
\For{$k=1\dots, K-1$}
\State $m_k = \min\{m, k\}$
\State $R_k=G(x_k)-x_k$
\State Solve
$\min_{\alpha^k=(\alpha^k_{k-m_k},\dots ,\alpha^k_{k})^\top}
\|\sum_{j=k-m_k}^{k}\alpha_j^k R_j\|_2^2$ subject to $\sum_{j=k-m_k}^{k} \alpha_j^k =1$
\State $x_{k+1}=(1-\beta)\sum_{j=k-m_k}^{k}\alpha_j^k x_j
        +\beta \sum_{j=k-m_k}^{k}\alpha_j^k G(x_j)$
\EndFor
\State \textbf{return} $x_{K}$
\end{algorithmic}
\end{algorithm}

One notable advantage of AA is that they do not require to know how the sequence is actually generated, thus the applicability of AA is quite wide. 
To apply AA, it is wise to keep in mind the following two aspects:  
\begin{itemize}
\item[(i)] The method only applies to the regime when the sequence admits a limit, which requires that the original iteration scheme has a convergence guarantee. 
\item[(ii)] The residual matrix $R$ can be rank-deficient, then instability may occur in computing $\alpha$. This would require the use of some regularization techniques to obtain a reliable $\alpha$ \cite{AA-optim3+}. By solving a regularized least square problem, $\alpha^k$ in Algorithm \ref{aagd} can be expressed as
\begin{equation}\label{ralpha}
[\alpha_{k-m_k}^k, \cdots, \alpha_{k-1}^k]= (U_k^\top U_k + \lambda I)^{-1}U_k^\top R_k,\quad\alpha^k_k = 1- \sum_{j=k-m_k}^{k-1} \alpha_j^k.
\end{equation}
where $U_k= [U_{k,k-m_k}, \cdots, U_{k,k-1}]$ with $U_{k,j} = R_k - R_j$, $\lambda>0$ is the regularization parameter introduced by adding 
$$
\lambda\|[\alpha^k_{k-m_k} \cdots \alpha^k_{k-1}]\|^2.  
$$
in line 5 of Algorithm \ref{aagd}. This treatment bears further comment, see the discussion on computing $\alpha$ in Section \ref{sec4.2}. 
\end{itemize}

\section{Anderson acceleration for optimization algorithms} 
The classical method to solve problem \eqref{probsetup} is Gradient Descent (GD), which is defined by
\begin{equation}\label{gd}
    x_{k+1} = x_k - \eta\nabla f(x_k),
\end{equation}
where $\eta>0$ is the step size. This can be viewed as the fixed-point iteration applied to $G(x)=x-\eta\nabla f(x)$. The fixed point $x^*$ of $G$, which corresponds to $\nabla f(x^*)=0$, is a stationary point of $f$. To illustrate how AA helps accelerate GD, we consider the quadratic minimization problem
\begin{equation}\label{quadf}
f(x)=\frac{1}{2}x^\top Ax-b^\top x,
\end{equation}
where $A$ is a symmetric positive definite matrix and the eigenvalues of A are bounded with $\sigma(A)\in[\mu,L]$. For this problem, GD reduces to the following fixed-point iteration: 
\begin{equation}\label{xke}
x_{k+1}-x^*=M(x_k-x^*),    
\end{equation}
where $M=I-\eta A$ and $x^*=A^{-1}b$. It can be verified that if $\eta=\frac{2}{L+\mu}$, then $\|M\|\leq\rho=\frac{L-\mu}{L+\mu}<1$, and $x_k$ converges to $x^*$ at a linear rate with
\begin{equation}\label{gd-cvg}
x_{i} - x^* = M^{i}(x_0 - x^*), \quad i=1, \cdots, k+1.   
\end{equation}
Now we investigate how the convergence rates can be improved using a linear combination of the previous iterations:
\begin{align*}
\|x_{k+1}^{AA}- x^*\| = \Bigg\|\sum_{i=0}^k \alpha_ix_{i} -x^*\Bigg\| & = \Bigg\|\sum_{i=0}^k \alpha_i (x_i - x^*)\Bigg\|= \| p_k(M)(x_0-x^*)\|,
\end{align*}
where we define the polynomial $p_k(M) = \sum_{i=0}^k\alpha_i M^i$ and $p_k(1)=1$. 
Using Chebyshev polynomials, the $l_2$ norm of the above error is bounded by (see Proposition 2.1 in \cite{AA-optim3+}): 
$$
\min_{p_k \in P^k, p_k(1)=1} \|p_k(M)(x_0-x^*)\|\leq \left\{ 
\begin{array}{ll}
\frac{2\gamma^k}{1+\gamma^{2k}}\|x_0-x^*\| & k<d,\\
0 & \text{otherwise}.
\end{array}
\right. 
$$
where $P^k$ is the subspace of polynomials of degree at most $k$ and   
$$
\gamma=\frac{1-\sqrt{1-\rho}}{1+\sqrt{1-\rho}}<1. 
$$
This when applied to the case $M=I-\eta A$ with $\eta=\frac{2}{L+\mu}$, we obtain 
$$
\gamma=\frac{\sqrt{\mu+L}-\sqrt{2\mu}}{\sqrt{\mu+L}+\sqrt{2\mu}}.
$$
This error bound is slightly smaller than that stated in \cite[Corollary 2.2]{AA-optim3+}. This is faster (since $\gamma<\rho$) than GD which based on (\ref{gd-cvg} admits the convergence rate of
$$
\|x_{k+1} - x^*\| = \rho^{k+1}\|x_0 - x^*\|.
$$
For AA with a finite $m$, how AA can help accelerate the convergence of an optimization method is relatively less understood. We present a simple result to show the accelerating effect.  
\begin{theorem}[Convergence of AA-GD for quadratic functions]\label{thm1}
For $f(x)=\frac{1}{2}x^\top Ax-b^\top x$ with  $A^\top=A$ and $\sigma(A)\in[\mu,L]$, Algorithm \ref{aagd} with $\beta=1$ and $\eta \leq \frac{2}{L+\mu}$ guarantees
\begin{equation}\label{fk-}
    \|\nabla f(x_{k+1})\|/\|\nabla f(x_{k})\|  \leq \delta_k (1-\eta\mu),
\end{equation}
where 
$\delta_k
= \|\Pi_k\nabla f(x_k)\|/\|\nabla f(x_k)\|
\leq 1$.
Here $\Pi_k$ is a projection operator defined by 
$
\Pi_k=I-U_k(U_k^\top U_k+\lambda I)^{-1}U_k^\top,
$
where $U_k= [U_{k,k-m_k}, \cdots, U_{k,k-1}]$ with $U_{k,j} =\nabla f(x_k)-\nabla f(x_{j})$, $\lambda$ is the regularization parameter for computing $\alpha$ (\ref{ralpha}). Moreover, we have the following convergence rate: 
\begin{equation}\label{fk} 
\|x_{k+1}^{AA}-x^*\| \leq \prod_{j=0}^{k} \delta_j (1-\eta\mu)^{k+1} \left(\frac{L}{\mu}\right)
\|x_0-x^*\|.
\end{equation}
\end{theorem}
The proof is deferred to Appendix \ref{pf1}. 
This result defines the gain of the optimization stage at iteration $k$ to be the ratio $\delta_k$ of the optimized gradient compared to that of the vanilla GD. 

The above result can be extended to the general non-quadratic case. 
\begin{theorem}[Convergence of AA-GD for non-quadratic functions]\label{thm2}
For $f(x)\in C^2$ and $\sigma(D^2f) \in [\mu, L]$, let $x^*$ be a minimum of $f(x)$, $x_k$ be the solution generated by Algorithm \ref{aagd} with $\beta=1$ and $\eta\leq\frac{2}{L+\mu}$. When $\|x_k-x^*\|$ is sufficiently small for a large
$k$, we have 
\begin{equation}\label{3.7} 
 \|\nabla f(x_{k+1})\|/\|\nabla f(x_k)\|
\leq \delta_k (1-\eta\mu)+o(k), 
\end{equation}
where $o(k)\to 0$ as $k\to \infty$,  with $\Pi_k$ defined in Theorem 3.1,  
\begin{equation*}
\delta_{k}=\|\Pi_k\nabla f(x_k)\|/\|\nabla f(x_k)\| \leq 1.
\end{equation*}
\end{theorem}
The proof is deferred to Appendix \ref{pf2}.
\begin{remark} \label{rem3.3} Note that the first term on the right hand side of (\ref{3.7}) is the same as that in (\ref{fk-}) for the quadratic case, and the second term converges to zero of order 
$$
\omega(\|x_{i}-x^*\|),
$$
where  $i=k-m_k, \cdots, k$, and $\omega(\cdot)$ is a modulus of continuity of $\nabla^2 f$. Under stronger assumptions that $\nabla^2f$ is also  Lipschitz continuous and $\eta=\frac{2}{L+\mu}$, the convergence bound given in  \cite[Theorem 2]{AA-optim4} corresponds to  
$$
 \|\nabla f(x_{k+1})\|/\|\nabla f(x_k)\| \leq 1-\eta \mu +O(\Delta_k), \quad  \Delta_k:= \max_{i\in [m_k]} \|x_k-x_{k-i}\|.
$$
\end{remark}

We note that the usual momentum acceleration differs from AA since the former seeks to modify the algorithmic structures at each iteration. For example, using two recent iterates, one has the following update rule: 
\begin{align*}
    &y_k = (1+\beta)x_k -\beta x_{k-1}, \\
    &x_{k+1}  = y_k - \eta \nabla f(y_k). 
\end{align*}
This when taking $\beta = (k-1)/(k+2)$ is the celebrated Nesterov's acceleration \cite{Nest1, Nest2}, which guarantees optimal convergence with sublinear rate $O(1/k^2)$ for convex functions $f(x)$ with L-Lipschitz-continuous gradients. With AA more historical iterates are used to improve the rate of convergence of fixed-point iterations.  
This said, there are also works on accelerating the Nesterov method by AA  \cite{AA-nest}.

\section{Anderson acceleration for AEGD}

Having seen the acceleration of AA for GD with the convergence gain quantified by a factor $\delta_k \leq 1$ at $k-$th step, we proceed to see how it may translate to improved asymptotic convergence behavior for AEGD. 

\subsection{Review of AEGD}
 Denote $F(x):= \sqrt{f(x)+c}$ with $c\in\R$ such that $f(x)+c>0$. AEGD is defined by 
\begin{subequations}\label{aegd}
\begin{align}
v_k &= \nabla F(x_k),\\
r_{k+1} &= \frac{r_k}{1+2 \eta v^2_k},\\
x_{k+1} &= x_k -2 \eta r_{k+1} v_k.
\end{align}
\end{subequations}
Here $v_k$ can be written explicitly as $\nabla f(x_k)/(2F(x_k))$, which represents the descent direction. 
This method is shown to be unconditional energy stable in the sense that $r_k$, serving as an approximation to $F(x_k)$, decreases monotonically regardless of the step size $\eta$. 

Note that (\ref{aegd}c) can be rewritten as 
\begin{equation}\label{aegd-2}
x_{k+1} = x_k - \eta_k \nabla f(x_k), \quad \eta_k = \eta \frac{r_{k+1}}{F(x_k)},
\end{equation} 
which shows that AEGD is not a fixed-point iteration due to the energy adaptation in the effective step size. However, since AA is a sequence acceleration technique, we expect AEGD to still gain a speedup from this technique as long as it is convergent. Indeed, it has been established in \cite{LT20} that under some structural conditions on $f$, the sequence $x_k$ generated by AEGD converges to a minimizer $x^*$. Below we present the convergence result of AEGD under the Polyak-Lojasiewicz (PL) condition. 

\begin{theorem}[Convergence of AEGD \cite{LT20}]
Suppose that the objective function $f$ is L-smooth: $\|\nabla^2 f(x)\|\leq L$ for any $x\in\R^n$ and $\max \|\nabla f(x)\|\leq G_\infty$. Denote $x^*\in \argmin_xf(x)$ and $r^*=\lim_{k\to\infty}r_k$, then
$$
r_k> r^* > F(x^*)(1-\eta/\tau),
$$
where $\tau>0$ is a constant that depends on $F(x^*), F(x_0), L$ and $G_\infty$. Moreover, for $f$ that further satisfies the PL condition: $\|\nabla f(x)\|^2\geq \mu (f(x)-f(x^*))$ for any $x\in\R^n$. If $\max_{j <k} \eta_j\leq 1/L$, then $x_k\to x^*$ and
\begin{align*}
    & \sum_{k=0}^\infty \|x_{k+1}-x_k\| \leq \frac{4}{\sqrt{2\mu}}\sqrt{f(x_{0})-f(x^*)}, \\
    & f(x_k)-f(x^*) \leq e^{-c_0kr_k} (f(x_{0})-f(x^*)), \quad c_0:=\frac{\mu\eta}{\sqrt{f(x_{0})+c}}.
\end{align*}
\end{theorem}

The convergence result stated above indicates that when $\eta$ is suitable small, we have 
$$
x_k \to x^*, \quad \eta_k=\eta \frac{r_{k+1}}{F(x_k)} \to \eta\frac{r^*}{F(x^*)}=\eta^* >0.  
$$
Hence asymptotically, AEGD is close to the fixed-point iteration with  $G(x)=x - \eta^* \nabla f(x)$ when the $x_k$ is close to $x^*$.




\subsection{Algorithm}\label{sec4.2}
Based on the above argument, we incorporate AA into AEGD and summarize the resulting algorithm, named AA-AEGD, in Algorithm \ref{aaegd}. 

\begin{algorithm}
\caption{AA-AEGD($m,q$) (proposed)} 
\label{aaegd}
\begin{algorithmic}[1] 
\Require $x_0$: initialization of the parameter, $c$: a constant satisfies $f(x)+c>0$; $r_0=\sqrt{f(x_0)+c}\bf{1}$: initialization of the energy parameter; $\eta$: step size, $m$: length of the window; $\beta$: relaxation parameter; $K$: total number of iterations.
\For{$k=0\dots, K-1$}
\State $v_k = \nabla f(x_k)/ (2\sqrt{f(x_k)+c})$
\State $r_{k+1} = r_k/ (1+2\eta v_k^2)$
\State $x_{k+1}=x_k-2\eta r_{k+1}v_k$
\State $R_k = x_{k+1}-x_k$
\If{$k=0 \;\text{mod}\; q$ and $k\neq0$}
\State Solve
$\min_{\alpha^k=(\alpha^k_{k-m},\dots ,\alpha^k_{k})^\top}
\|\sum_{j=k-m}^{k}\alpha_j^k R_j\|^2_2$ subject to $\sum_{j=k-m}^{k} \alpha_j^k =1$
\State $x_{k+1}^{AA}=(1-\beta)\sum_{j=k-m}^{k}\alpha_j^k x_j
        +\beta \sum_{j=k-m}^{k}\alpha_j^k x_{j+1}$
\State $x_{k+1}=x_{k+1}^{AA}$
\EndIf
\EndFor
\State \textbf{return} $x_K$
\end{algorithmic}
\end{algorithm}

Some implementation techniques in the algorithm will be used in our experiments. We thus introduce and discuss their effects through experiments on a quadratic function.

{\bf Computation of $\alpha$.}
The optimization problem in Algorithm \ref{aaegd} can be formulated as an unconstrained least-squares problem: 
\begin{equation}\label{ucls}
    \min_{(\alpha^k_{k-m},\dots ,\alpha^k_{k-1})^\top}
    \Bigg\|R_k-\sum_{j=k-m}^{k-1}\alpha_j^k (R_k-R_j)\Bigg\|^2_2,\quad\text{set}\quad \alpha_k^k = 1-\sum_{j=k-m}^{k-1}\alpha_j^k,
\end{equation}
which can be solved by standard methods such as QR decomposition with a cost of $\mathcal{O}(2m^2n)$, where $n$ is the dimension in state space. Since $m$ is typically used less than $10$ in practice, AA can be implemented with efficient computation. In our experiments, we add a regularization of $\lambda\|[\alpha^k_{k-m} \cdots \alpha^k_{k-1}]\|^2_2$ to \eqref{ucls} to avoid singularity, where $\lambda > 0$ is a parameter, as implemented in \cite{AA-optim3+}. The explicit expression of $\alpha^k$ is given in (\ref{ralpha}).

\begin{figure}[ht]
\begin{subfigure}[b]{.5\textwidth}
  \centering
  \includegraphics[width=1\linewidth]{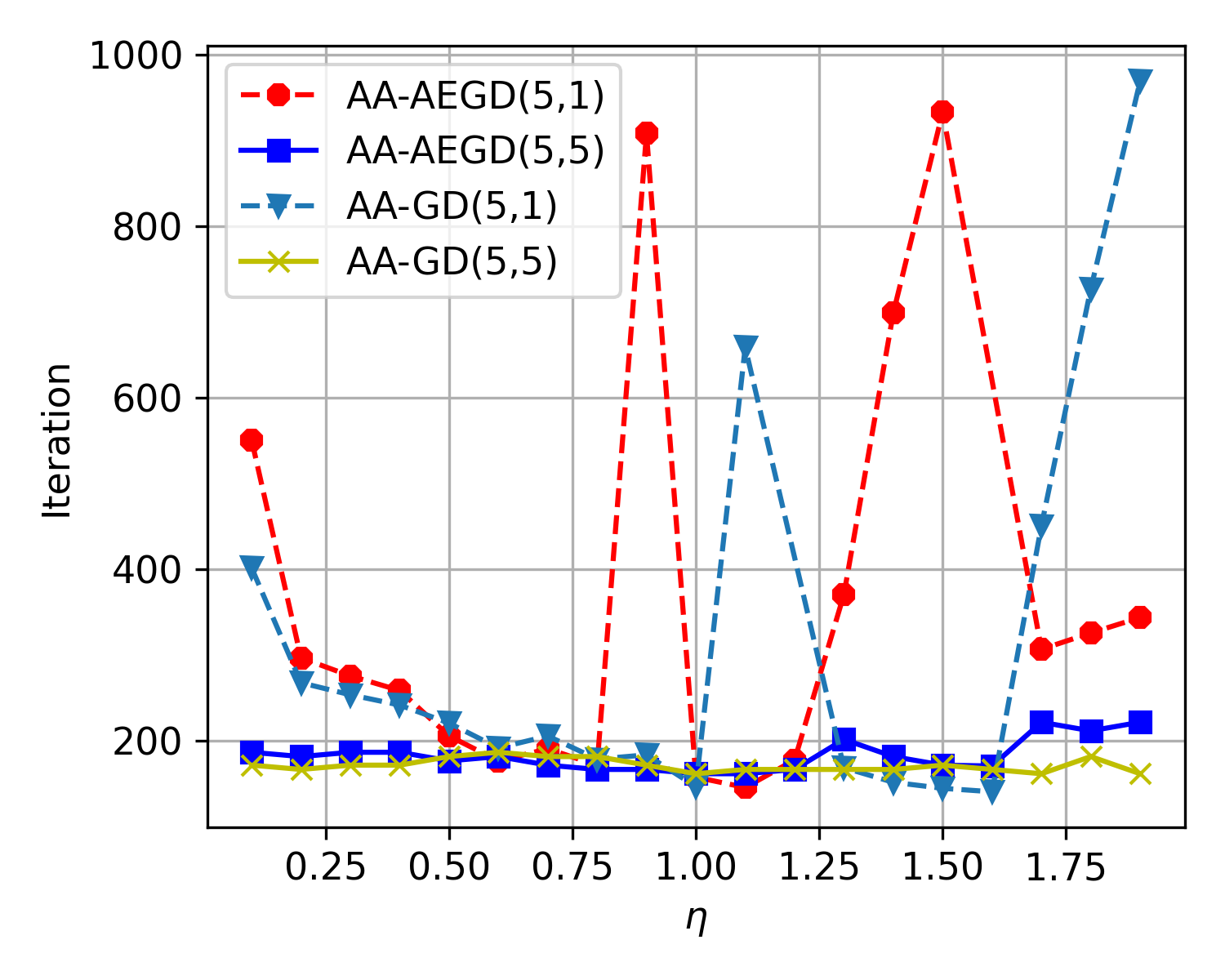}
  \caption{The effect of $\eta$}
\end{subfigure}%
\begin{subfigure}[b]{.5\textwidth}
  \centering
  \includegraphics[width=1\linewidth]{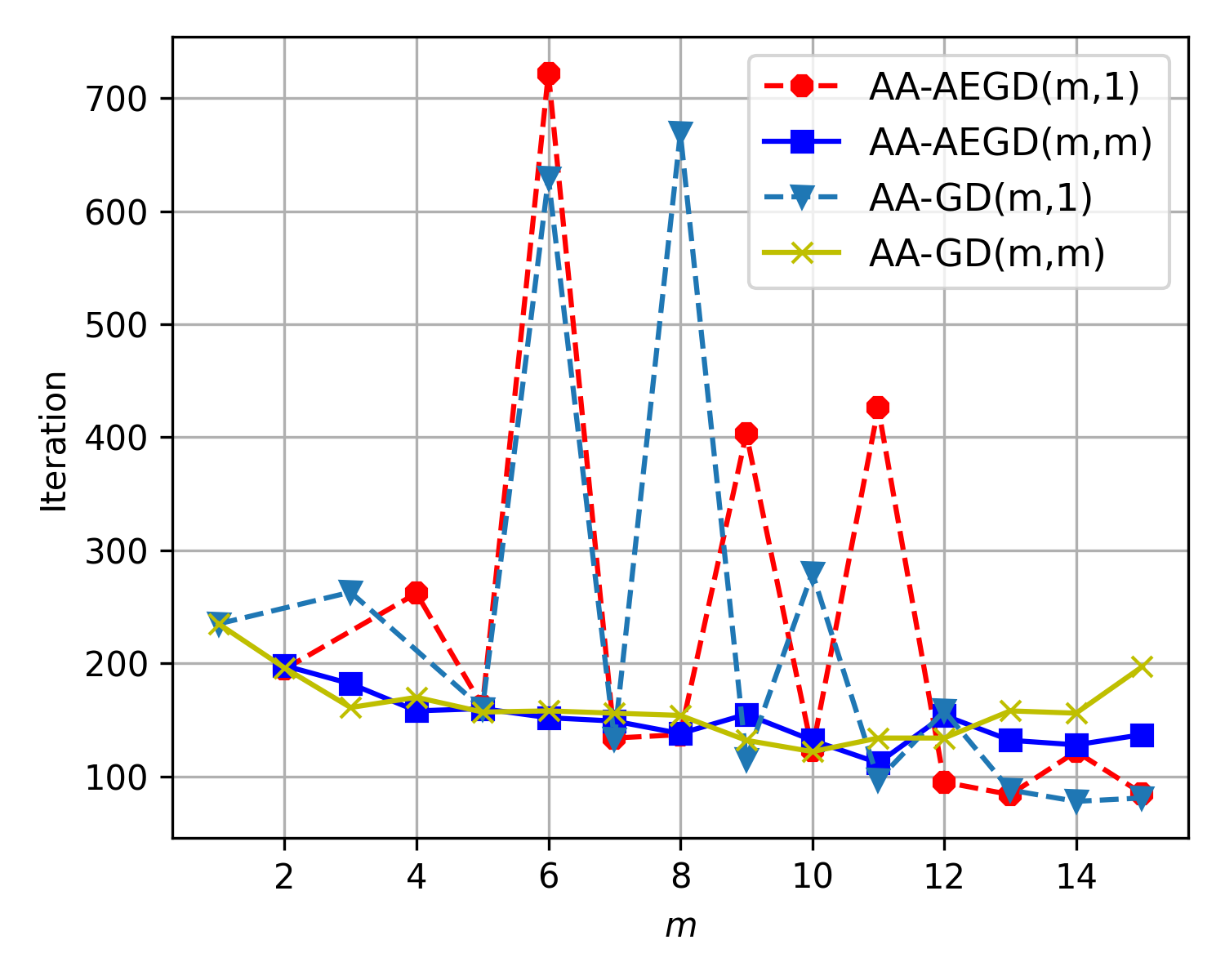}
  \caption{The effect of $m$}
\end{subfigure}
\caption{The effect of $\eta$ and $m$ on AA-GD/AEGD($m,1$) and AA-GD/AEGD($m,q$) for solving 
the quadratic problem (\ref{quadf}) with $\kappa=10^3$, where $\kappa=\sigma_{max}(A)/\sigma_{min}(A)$ is the condition number}
\label{fig-quad_eta_m}
\end{figure}

{\bf AA implementation of frequency.} 
Unlike the standard version presented in Algorithm \ref{aagd}, which implements AA in every step, Algorithm \ref{aaegd} implements AA in every $q$ steps. We denote the former version as AA-GD/AEGD($m,1$), and the latter version as AA-GD/AEGD($m,q$). Obviously, AA($m,q$) requires less computation than AA($m,1$). Another advantage of AA($m,q$) is that it is less sensitive to the effect of $\eta$ and $m$, as evidenced by  numerical experiments (see Figure \ref{fig-quad_eta_m}). This is a favorable property since it indicates that the proposed algorithm requires little tuning of hyperparameters.



\subsection{Application to proximal gradient algorithm}

Consider the following minimization problem 
\begin{equation}\label{cpf}
\min_{x\in\R^n} L(x) = f(x) + h(x),    
\end{equation}
where $f:\R^n\to\R$ is $L$-Lipschitz smooth: $\|\nabla f(x)-\nabla f(y)\|_2\leq L\|x-y\|_2$, and $h$ is a proper closed and convex function. In many applications, such as data analysis and machine learning, $f(x)$ is a convex data fidelity term, and $h(x)$ is a certain regularization.  For example, we take $h(x)=\frac{1}{2}\|x\|^2_2$ in Section \ref{logistic} and Section \ref{lasso}. A classical method for solving (\ref{cpf}) is the Proximal Gradient Algorithm (PGA):
\begin{subequations}\label{pga}
\begin{align}
y_{k+1} &= x_k - \eta \nabla f(x_k),\\
x_{k+1} &= \text{Prox}_{\eta h}(y_{k+1}),
\end{align}    
\end{subequations}
where the proximal operator associated with $h$ is defined as 
\begin{equation}
\text{Prox}_{\eta h}(y) = \argmin_{x}\Bigg\{h(x)+\frac{1}{2\eta}\|x-y\|^2_2\Bigg\}.
\end{equation}
To apply AA to PGA, \cite{AA-PGD} treats (\ref{pga}a) as the fixed-point iteration $y_{k+1}=G(y_k)$ with 
\begin{equation}
G(y) = \text{Prox}_{\eta h}(y) - \eta\nabla f(\text{Prox}_{\eta h}(y)).
\end{equation}
Here AA is applied on the auxiliary sequence $\{y_k\}$, and the authors in  \cite{AA-PGD}  show that  
fast convergence of $\{y_k\}$ will lead to fast convergence of $\{x_k\}$.  

We apply AA-AEGD to PGA 
to take advantage of both AA and energy adaptation on step size. We summarize the resulting scheme in Algorithm \ref{aapeg}. Following  \cite{AA-PGD}, we add a step checking (Line 11) to ensure convergence. The numerical results on the performance of Algorithm \ref{aapeg} are presented in Section \ref{logistic} and Section \ref{lasso}.

\begin{algorithm}
\caption{AA- AEGD($m,q$) for solving problem (\ref{cpf}).} 
\label{aapeg}
\begin{algorithmic}[1] 
\Require $x_0$: initialization of the parameter; $y_0=x_0$;
$c$: a constant satisfies $f(x)+c>0$; 
$r_0=\sqrt{f(x_0)+c}\bf{1}$: initialization of the energy parameter; 
$\eta$: step size;
$m$: length of the window; 
$\beta$: relaxation parameter; 
$K$: total number of iterations.
\For{$k=0\dots, K-1$}
\State $v_k = \nabla f(x_k)/ (2\sqrt{f(x_k)+c})$
\State $r_{k+1} = r_k/ (1+2\eta v_k^2)$
\State $y_{k+1}=x_k-2\eta r_{k+1}v_k$
\State $x_{k+1}=\text{Prox}_{\eta h}(y_{k+1})$
\State $R_k = y_{k+1}-y_k$ 
\If{$k=0 \;\text{mod}\; q$ and $k\neq0$}
\State Solve
$\min_{\alpha^k=(\alpha^k_{k-m},\dots ,\alpha^k_{k})^\top}
\|\sum_{j=k-m}^{k}\alpha_j^k R_j\|^2_2$ subject to $\sum_{j=k-m}^{k} \alpha_j^k =1$
\State $y_{k+1}^{AA}=(1-\beta)\sum_{j=k-m}^{k}\alpha_j^k y_j
        +\beta \sum_{j=k-m}^{k}\alpha_j^k y_{j+1}$
\State $x_{k+1}^{AA}=\text{Prox}_{\eta h}(y_{k+1}^{AA})$
\If{$f(x_{k+1}^{AA}) \leq f(x_k) - \frac{\eta}{2}\|\nabla f(x_k)\|^2_2$}
\State $x_{k+1}=x_{k+1}^{AA}$
\EndIf
\EndIf
\EndFor
\State \textbf{return} $x_K$
\end{algorithmic}
\end{algorithm}

\section{Numerical Experiments}\label{exp}

In this section, we demonstrate the performance of AA-AEGD on several optimization problems, including the 2D Rosenbrock function and two machine learning problems. For parameters, we set $c=1$ for AEGD and AA-AEGD, $\beta=1$ and $q=m$ for AA-GD/AEGD($m,q$). All experiments are conducted in Python and run on a laptop with four 2.4 GHz cores and 16 GB of RAM. \footnote{The code is available at \url{https://github.com/howardjhe/AA-AEGD.git}.}

In Section \ref{nc}, the effect of $\eta$ and $m$ on AA-GD/AEGD($m,q$) is studied in the non-convex case. Then the performance of AA-AEGD (Algorithm \ref{aaegd}) is compared with GD, AEGD, and AA-GD. In Section \ref{logistic} and Section \ref{lasso}, we further implement AA-AEGD (Algorithm \ref{aapeg}) to solve problems in the form of (\ref{cpf}). We employ the large-scale and ill-conditioned datasets, Madelon ($\kappa \approx 2.1\times 10^4$) and MARTI0 ($\kappa \approx 7.5 \times 10^3$), obtained from the UCI Machine Learning Repository\footnote{The Madelon dataset is available from \url{https://archive.ics.uci.edu/ml/datasets/}} and ChaLearn Repository\footnote{The MARTI0 dataset is available from \url{http://www.causality.inf.ethz.ch/home.php}}, respectively. As for the parameters used in the competitors, including PGA, APGA, and AA-PGA, we keep all the same settings as in \cite{AA-PGD}.

\subsection{Non-convex function} \label{nc}
\begin{figure}[ht]
\begin{subfigure}[b]{.5\textwidth}
  \centering
  \includegraphics[width=1\linewidth]{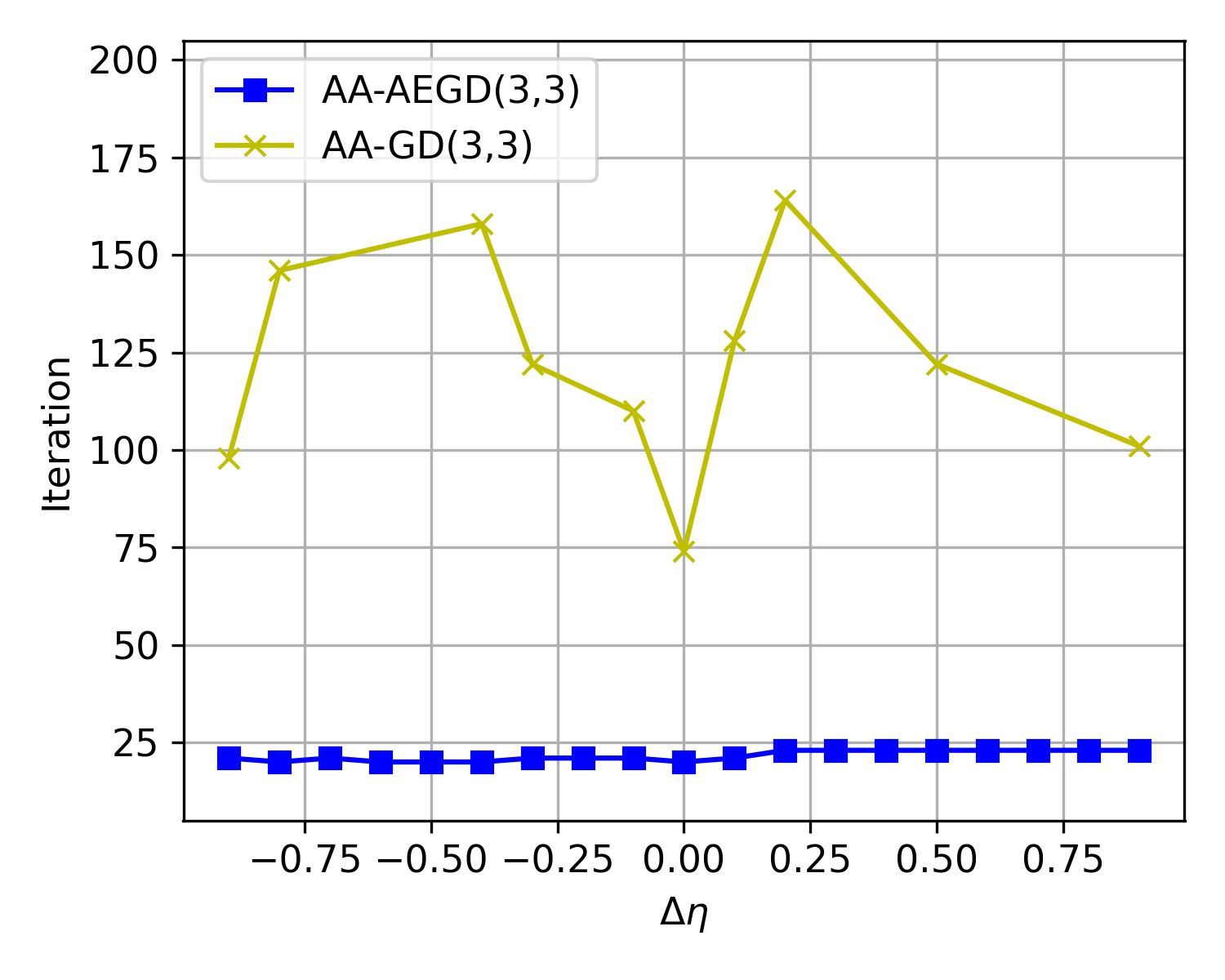}
  \caption{The effect of $\eta$}
\end{subfigure}%
\begin{subfigure}[b]{.5\textwidth}
  \centering
  \includegraphics[width=1\linewidth]{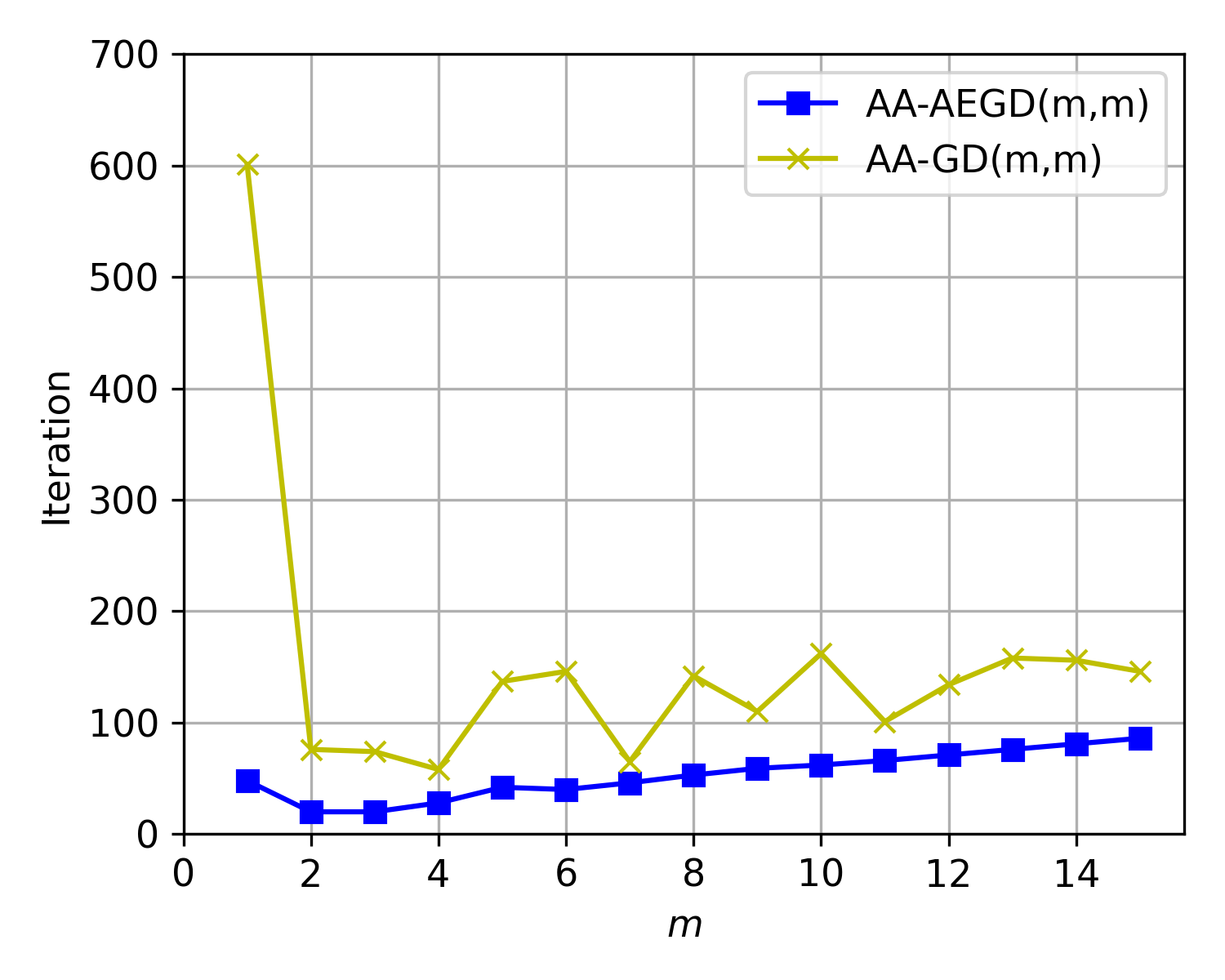}
  \caption{The effect of $m$}
\end{subfigure}
\begin{subfigure}[b]{.5\textwidth}
  \centering
  \includegraphics[width=1\linewidth]{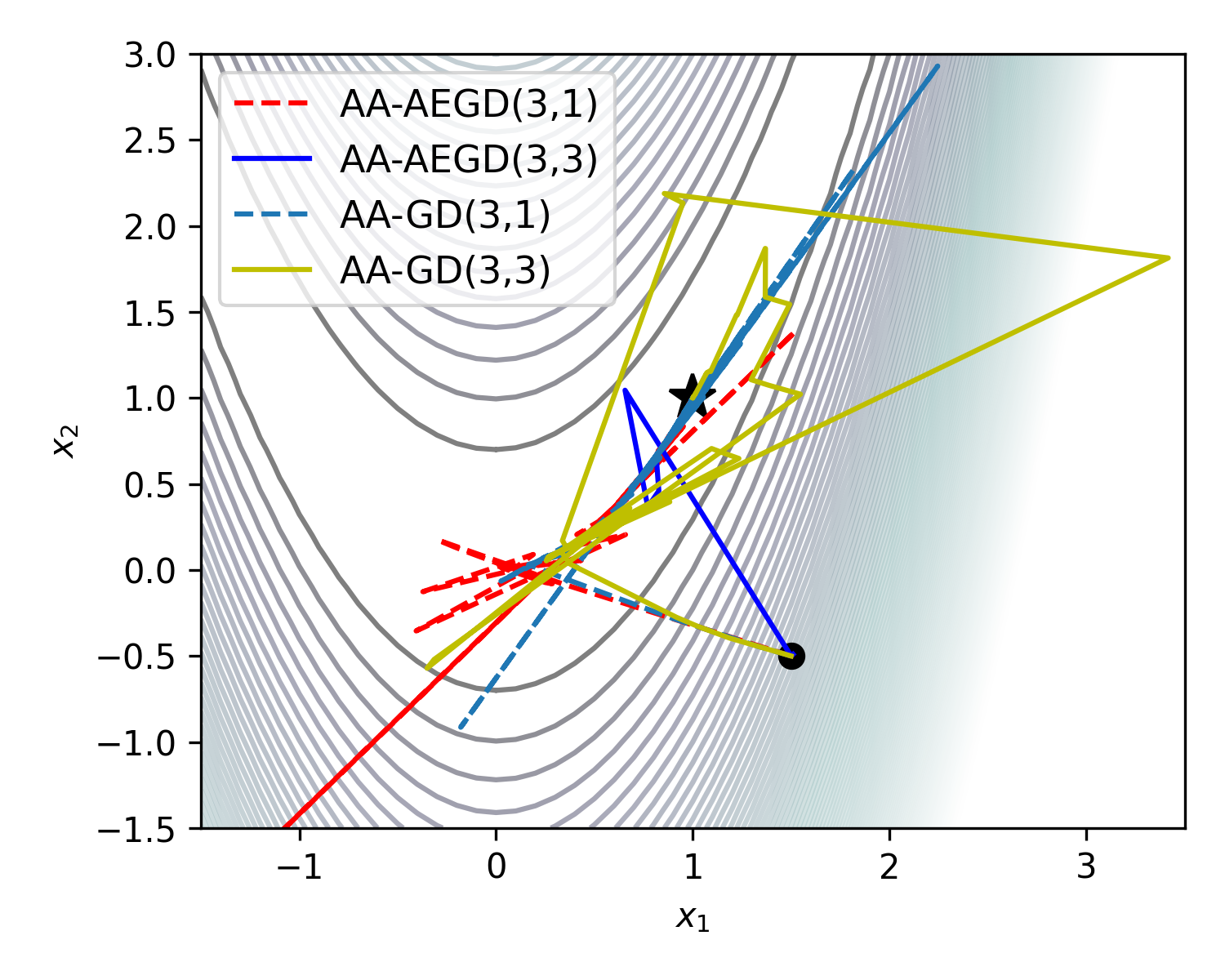}
  \caption{The effect of frequency implementation}
\end{subfigure}%
\begin{subfigure}[b]{.5\textwidth}
  \centering
  \includegraphics[width=1\linewidth]{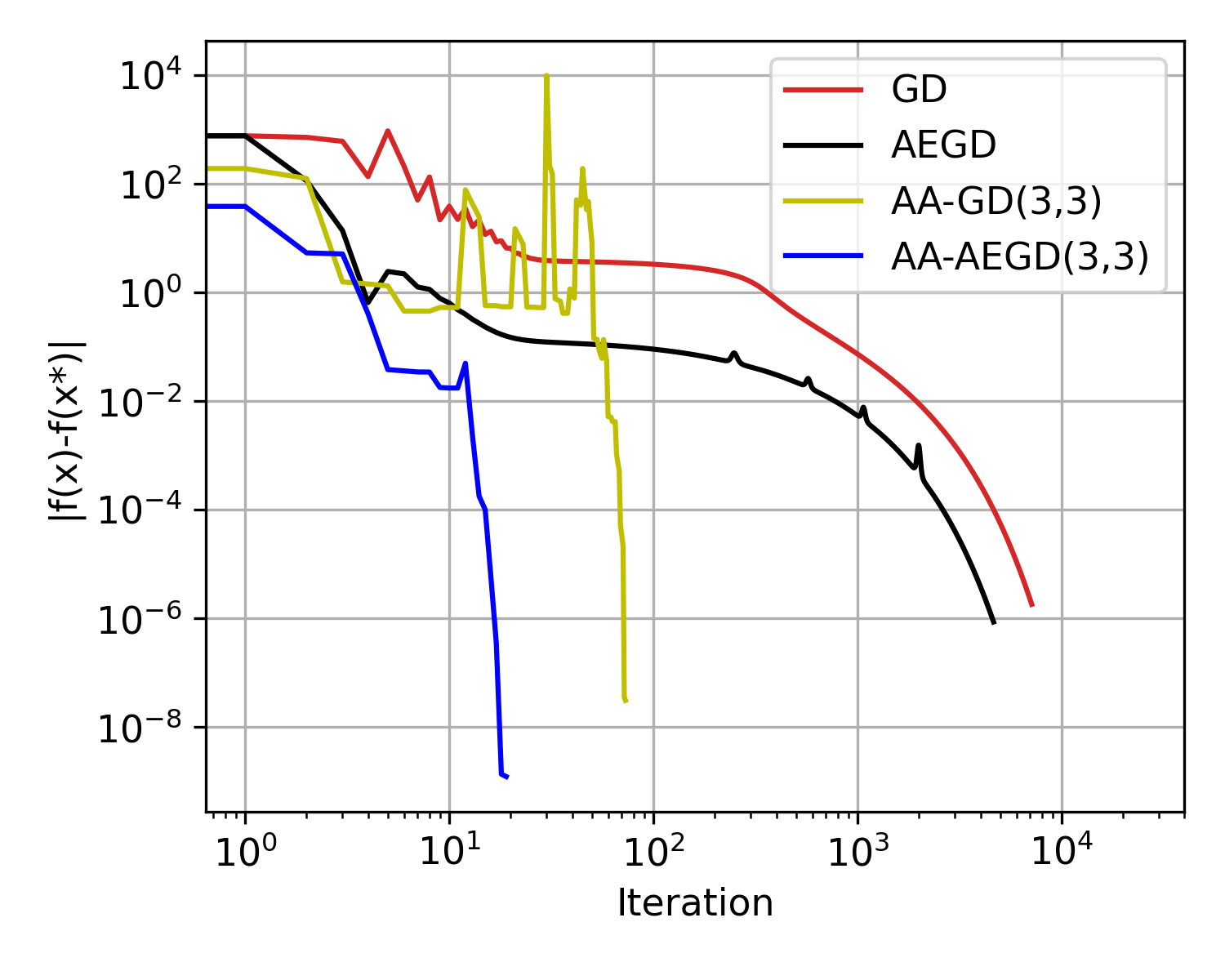}
  \caption{The number of iterates}
\end{subfigure}
\caption{Comparison of AA-AEGD and AA-GD on the 2D Rosenbrock function. The initial point is set at $(1.5,-0.5)$.}
\label{fig-rosen_eta_m}
\end{figure}

Consider the following 2D Rosenbrock function
\begin{equation*}
f(x_1,x_2)=(1-x_1)^2+100(x_2-x_1^2)^2.
\end{equation*}
This is a non-convex function with a global minimum $(1,1)$.  
As shown in Figure \ref{fig-rosen_eta_m} (c), the global minimum (marked as the black star) is inside a long, narrow valley. It takes only several iterations to find the valley, while a lot more iterations are needed to converge to the global minimum.

We first study the effect of $\eta$ and $m$ on AA-GD/AEGD($m,q$). Since the feasible ranges of $\eta$ for GD and AEGD are different, we take different intervals of $\eta$ for the two methods when studying the effect of $\eta$ on AA-GD/AEGD($m,q$). The interval of $\eta$ for each method is determined as follows: we first search for $\bar\eta$ with which the method reaches the minimizer with the least iterations; then take $[\bar\eta-\Delta\eta, \bar\eta+\Delta\eta]$ as the study interval of $\eta$. In our experiment, we find $\bar\eta_1 = 1.9\times 10^{-4}$ for GD, and $\bar\eta_2 = 6.4\times 10^{-3}$ for AEGD. For both methods, we set $\Delta \eta = 1.9$. The results are presented in Figure \ref{fig-rosen_eta_m} (a) (b). 




Figure \ref{fig-rosen_eta_m} (c) presents a comparison of trajectories of AA-GD/AEGD($m,1$) and AA-GD/AEGD($m,q$). We observe that compared with AA-GD/AEGD(3,3), the paths of AA-GD/AEGD(3,1) detour far away from the minimizer. Also, AA-AEGD(3,3) reaches the minimizer with the shortest path compared with others. A comparison of GD, AEGD, and AA-GD, AA-AEGD is given in Figure \ref{fig-rosen_eta_m} (d).

\subsection{Constrained logistic regression problem}\label{logistic} Consider the logistic regression problem with a bounded constraint:
\begin{align*}
&\min_{x\in \mathbb{R}^n} \frac{1}{M}\sum_{i=1}^M \log(1+\exp(-y_i a_i^{\top}x)) + \mu\|x\|_2^2,\\ 
&\mbox{subject to}\ \|x\|_\infty \leq 1,
\end{align*}
where $a_i\in \mathbb{R}$ are sample data and $y_i\in \{-1,1\}$ are corresponding labels to the samples. We set $\mu=10$ and apply $\lambda = 10^{-10}$ as the regularization parameter to \eqref{ucls}. For this  problem, we set the initial point at $\bm{0}$, and keep the step size, $\eta = 1/L_1$, where $L_1 = \|A\|_2^2/4M$ with $A=[a_1,\cdots,a_M]$, for PGA, APGA, and AA-PGA, and $m=5$ for AA-PGA as in \cite{AA-PGD}. We apply $\eta = 2/L_1$ to AEGD on both the Madelon and MARTI0 datasets. $\eta = 3/L_1$ and $4/L_1$ are adopted in AA-AEGD (Algorithm \ref{aapeg}) on the Madelon and MARTI0 datasets, respectively.

\begin{figure}[ht]
\begin{subfigure}{.5\textwidth}
  \centering
  \includegraphics[width=1.\linewidth]{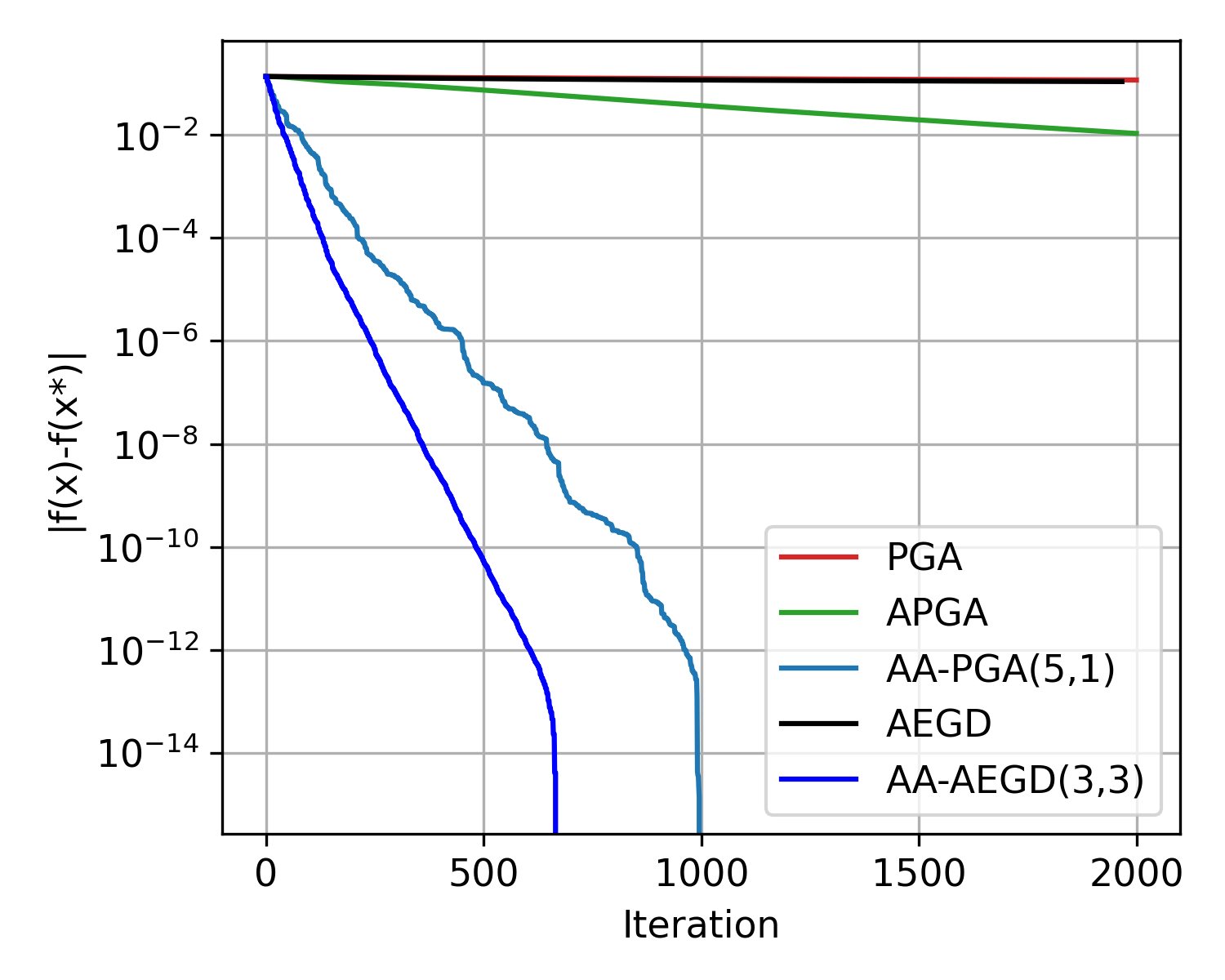}
  \caption{Madelon}
\end{subfigure}%
\begin{subfigure}{.5\textwidth}
  \centering
  \includegraphics[width=1.\linewidth]{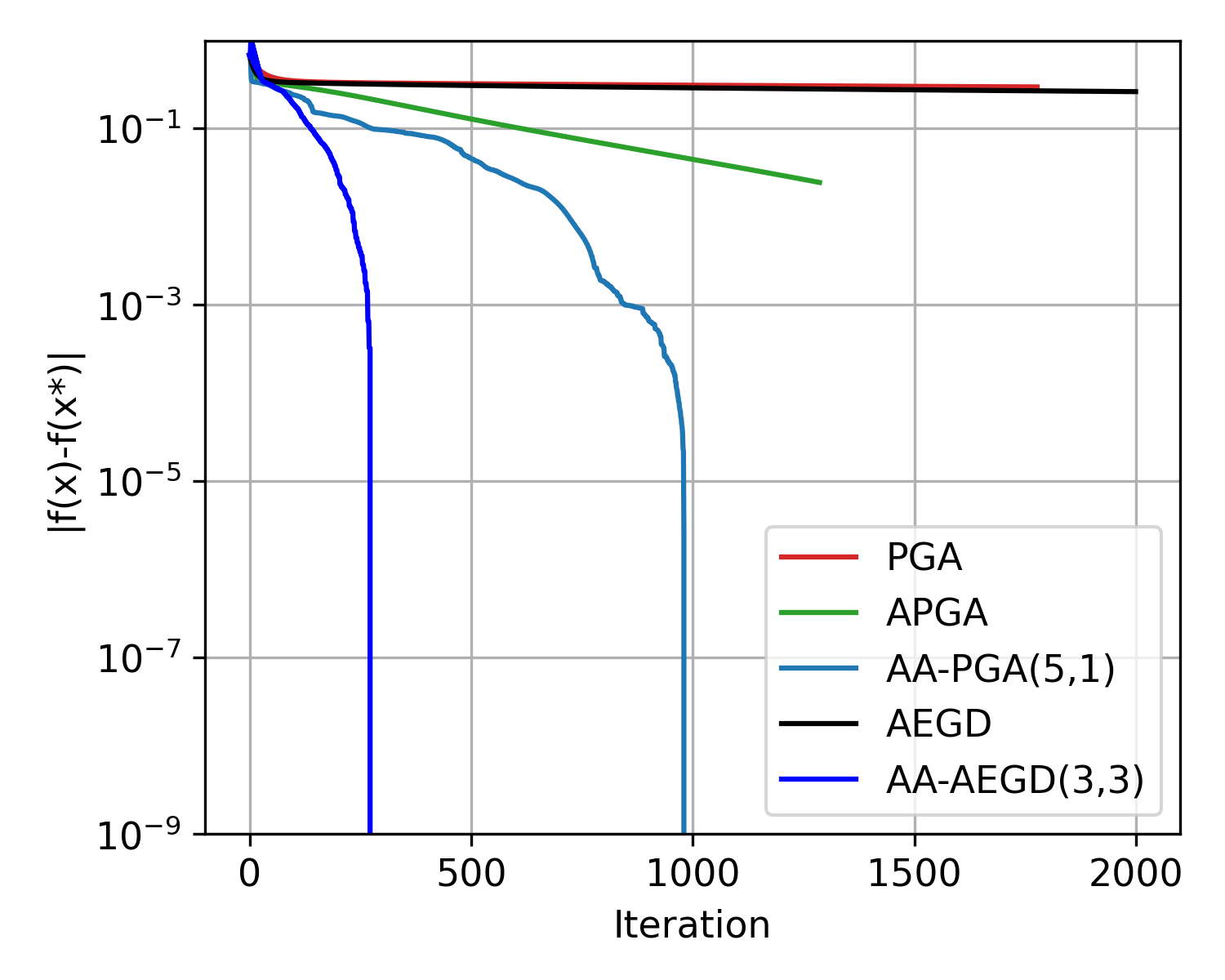}
  \caption{MARTI0}
\end{subfigure}
\caption{Constrained logistic regression problem}
\label{fig-log}
\end{figure}

Figure \ref{fig-log} demonstrates the superior performance of AA-AEGD on the constrained logistic regression problem compared with PGA, APGA, and AA-PGA. 

\subsection{Nonnegative least squares}\label{lasso} Consider the nonnegative least squares problem
\begin{align*}
& \min_{x\in \R^n} \frac{1}{2M}\|Ax-b\|^2_2 + \mu\|x\|^2_2\\
& \mbox{subject to}\ x \geq 0,   
\end{align*}
where $A\in \mathbb{R}^{M \times n}$ is the sample data and  $b\in \mathbb{R}^M$ are corresponding labels to the samples. We take $\mu=0.1$  and the initial point as $\bf{0}$. In this problem, we keep the step size, $\eta = 1/L_2$, where $L_2 = \|A\|_2^2/M$, for PGA, APGA, and AA-PGA, and $m=5$ for AA-PGA as same in \cite{AA-PGD}. We apply $\eta = 9/L_2$ to AEGD for both datasets. $\eta = 9/L_1$ and $3/L_1$ are employed with AA-AEGD on the Madelon and MARTI0 datasets, respectively. Overall, AA-AEGD exhibits faster convergence than AA-GD.


\begin{figure}[ht]
\begin{subfigure}{.5\textwidth}
  \centering
  \includegraphics[width=1.\linewidth]{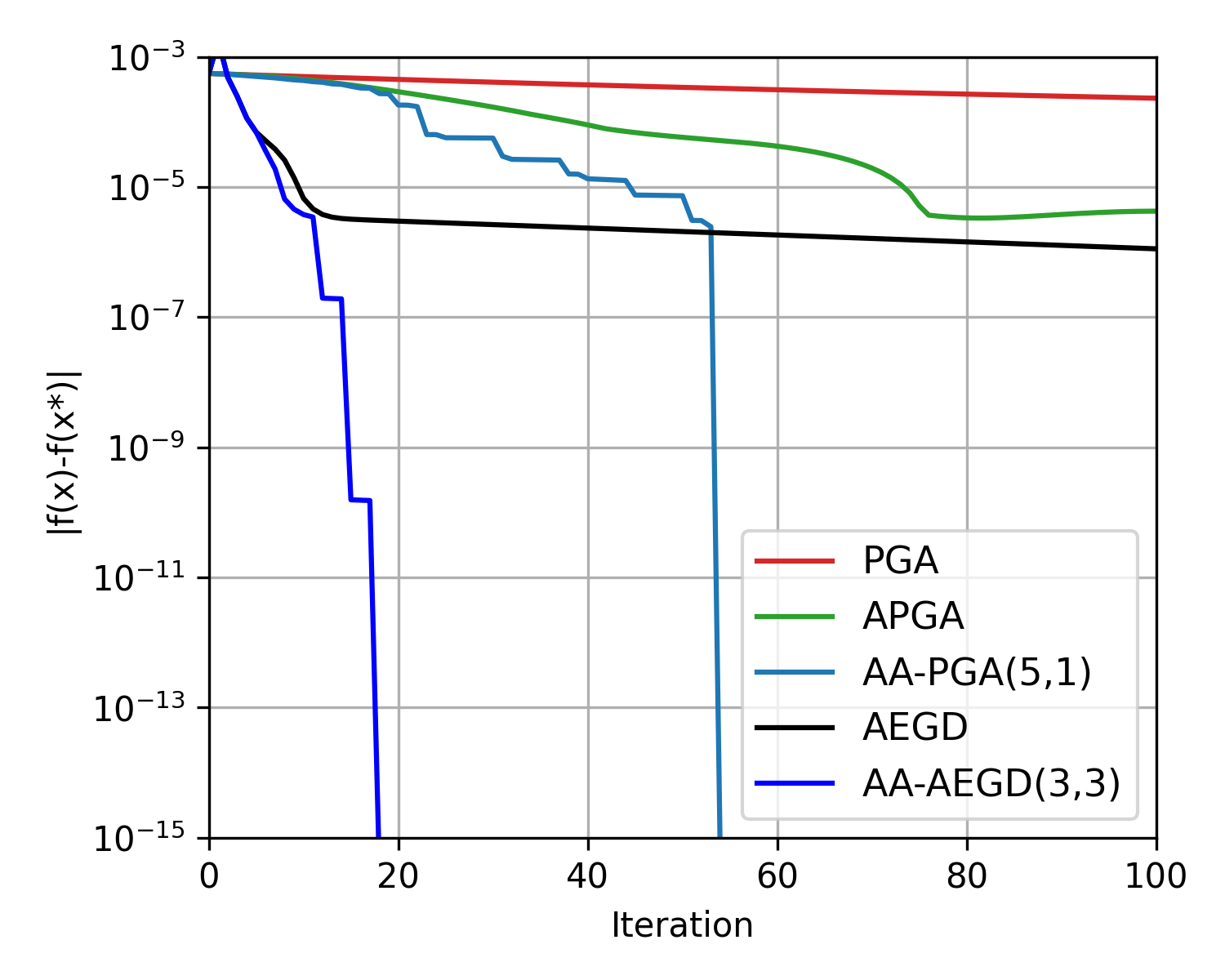}
  \caption{Madelon}
\end{subfigure}%
\begin{subfigure}{.5\textwidth}
  \centering
  \includegraphics[width=1.\linewidth]{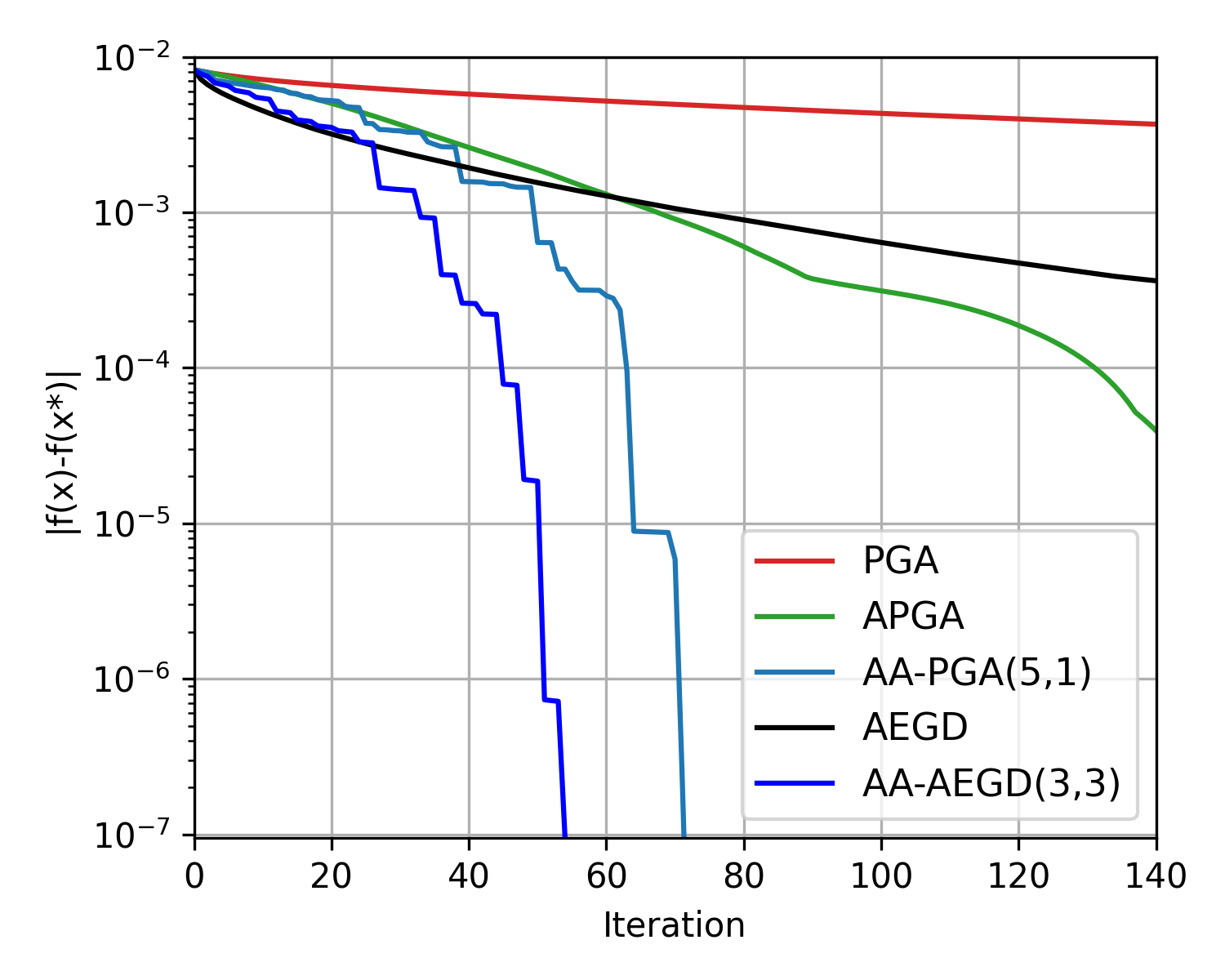}
  \caption{MARTI0}
\end{subfigure}
\caption{Nonnegative least-squares problem}
\label{fig-las}
\end{figure}

\section{Conclusion}
In this paper, we developed a novel algorithm of using Anderson acceleration for the adaptive energy gradient descent (AEGD) in solving optimization problems. First, we analyzed the gain in convergence rate of AA-GD for quadratic and non-quadratic optimization problems, and then we explained why AA can still be exploited with AEGD though it is not a fixed-point iteration. 
We deployed AA-AEGD for solving convex and non-convex problems and obtained superior performance on a set of  test problems. Moreover, we applied our algorithm to the proximal gradient method on two machine learning problems and observed improved convergence as well. Experiments show encouraging results of our algorithm and confirm the suitability of Anderson acceleration for accelerating AEGD in solving optimization problems.  With the easy implementation and excellent performance of AA-AEGD, it merits further study with more refined convergence analysis.  

\appendix
\renewcommand{\theequation}{\thesection.\arabic{equation}}

\section{Proof of Theorem \ref{thm1}} \label{pf1}
Using $G(x)=x-\eta\nabla f(x)$ and $\nabla f(x)= Ax-b$, we have
\begin{align*}
    \nabla f(x_{k+1}) &= Ax_{k+1}-b \\
    &= A\left(\sum\limits_{j=k-m}^k
    \alpha_j^kG(x_j)\right)-b\\
    &= A \left(\sum\limits_{j=k-m}^k
    \alpha_j^k(x_j-\eta(Ax_j - b))\right)-b \\
    &= \sum\limits_{j=k-m}^k\alpha_j^k [Ax_j - 
    \eta A(Ax_j-b)-b ]\\
    &= \sum\limits_{j=k-m}^k\alpha_j^k (I-\eta A)(Ax_j-b)\\
    &= (I-\eta A)\sum\limits_{j=k-m}^k\alpha_j^k\nabla f(x_j).
\end{align*}

Taking the $l^2$-norm on both sides, we get
\begin{equation}\label{gk}
    \|\nabla f(x_{k+1})\| \leq \|I-\eta A\| 
    \left\|\sum\limits_{j=k-m}^k \alpha_j^k \nabla f(x_j)\right\|. 
\end{equation}
Thus by taking $\eta \leq \frac{2}{L+\mu}$, (\ref{gk}) with the Anderson selection of $\alpha_j^k$ becomes
\begin{equation}\label{LuSk}
    \|\nabla f(x_{k+1})\| \leq  (1-\eta\mu)  \min_{\alpha^k } \left\|\sum\limits_{j=k-m}^k \alpha_j^k \nabla f(x_j)\right\|.
\end{equation}
Note that $\sum_{j=k-m}^k \alpha_j^k=1$ we have 
\begin{equation}\label{lsp}
\begin{aligned}
    \left\|\sum_{j=k-m}^k\alpha_j^k \nabla f(x_j)\right\|
    &= \left\|\nabla f(x_k)-\sum_{j=k-m}^{k-1} \alpha_j^k(\nabla f(x_k)-\nabla f(x_j))\right\| \\
    &= \left\|\nabla f(x_k)-U_k\alpha^k\right\|, \quad 
    \alpha^k:=[\alpha_{k-m}^k, \cdots, \alpha_{k-1}^k].
\end{aligned}    
\end{equation}
By solving the regularized least square problem:
\begin{equation}
\min_{\alpha \in \mathbb{R}^m} \frac{1}{2} \left\|\nabla f(x_k)-U_k\alpha\right\|^2 + \frac{\lambda}{2}\|\alpha\|^2,
\end{equation}
we obtain an explicit expression for the optimal weight vector:
\begin{equation}\label{lsalpha}
\alpha^k_*= (U_k^\top U_k+\lambda I)^{-1}U_k^\top\nabla f(x_k),\quad\alpha^k_k = 1- \alpha^k_*\cdot {\bf 1}.
\end{equation}
Hence 
\begin{equation*}
\left\|\sum_{j=k-m}^k\alpha^k_j \nabla f(x_j)\right\|  = \|\Pi_k\nabla f(x_k)\|. 
\end{equation*}
We thus have 
\begin{equation*}
    \|\nabla f(x_{k+1})\| \leq  (1-\eta\mu)  \|\Pi_k\nabla f(x_k)\| \leq \delta_k (1-\eta\mu)\|\nabla f(x_{k})\|.
\end{equation*}
This estimate can be further used to derive the bound on $\|x_{k+1}-x^*\|$. 
Note that 
$
\nabla f(x_{k+1}) = Ax_{k+1}-b = A(x_{k+1}-x^*), 
$
we then have
\begin{align*}
    \|x_{k+1}-x^*\| 
    &= \|A^{-1}(A(x_{k+1}-x^*))\| \\
    &\leq \|A^{-1}\| \|\nabla f(x_{k+1})\| \\
    &\leq \|A^{-1}\| C_k(1-\eta\mu)^{k+1} \|A(x_0-x^*)\| \\
    &\leq C_k(1-\eta\mu)^{k+1} \|A\|\|A^{-1}\| \|(x_0-x^*)\|,
\end{align*}
where $\|A\|\|A^{-1}\|$ is the condition number, and 
$C_k=\prod_{j=0}^{k} \delta_j$.  Thus, we have
$$
\|x_{k+1}-x^*\| \leq C_k(1-\eta\mu)^{k+1} \left(\frac{L}{\mu}\right)
\|(x_0-x^*)\|. 
$$
The proof is complete.

\section{Proof of Theorem \ref{thm2}}\label{pf2}

Denote $g(\cdot)=\nabla f(\cdot)$ and $g_k=\nabla f(x_k)$, the update of AA-GD is given by
\begin{equation*}
    x_{k+1}=\sum_{j=k-m}^{k} \alpha^k_j(x_j-\eta g_j). 
\end{equation*}
Set $y_k:= \sum_{j=k-m}^{k} \alpha^k_j x_j$, we make the following split: 
\begin{equation}\label{gab}
\eta g_{k+1} 
= a_k+b_k,
\end{equation}
where
\begin{align*}
a_k &= \eta g_{k+1}-x_{k+1}-(\eta g(y_k)-y_k),\\
b_k &= x_{k+1}+(\eta g(y_k)-y_k).
\end{align*}
We proceed to bound the $l^2$-norm of $a_k$ and $b_k$. First we have 
\begin{equation}\label{bda}
    \begin{aligned}
        \|a_k\| &= \|\eta(g_{k+1}-g(y_k))-(x_{k+1}-y_k)\| \\
        &= \|(\eta D^2f(\cdot)-I)(x_{k+1}-y_k)\| \quad
        \mbox{where}\ D^2f(\cdot)=\int_0^1 D^2f(tx_{k+1}+(1-t)y_k)dt\\
        &\leq \|\eta D^2f(\cdot)-I\| \|x_{k+1}-y_k\| \\
        &= \eta\|\eta D^2f(\cdot)-I\| \left\|\sum_{j=k-m}^{k} \alpha^k_j(-\eta g_j)\right\| \\
        &= \eta\|\eta D^2f(\cdot)-I\|\delta_k\|g_k\| \quad
        \mbox{where}\ \delta_k:=\frac{\|\Pi_k g_k\|}{\|g_k\|}\\
       &\leq \eta(1-\eta\mu)\delta_k\|g_k\|.
    \end{aligned}
\end{equation}
The last inequality is guaranteed by setting $\eta\leq\frac{2}{\mu+L}$. Recall that $\alpha^k$ is chosen by minimizing (\ref{lsp}), which ensures $\delta_k=\|\Pi_k g_k\|/\|g_k\| \leq 1$.

Using $\alpha^k$ defined in (\ref{lsalpha}), we now bound  $b_k$ as follows. 
\begin{equation*}
    \begin{split}
        b_k &= \sum_{j=k-m}^{k} \alpha^k_j(x_j-\eta g_j)
        -y_k+\eta g(y_k)\\
        &= x_k-y_k-\eta(g_k-g(y_k)) +\sum_{j=k-m}^{k-1} \alpha^k_j((x_j-\eta g_j)-(x_k-\eta g_k))\\
        &= (I-\eta D^2f(\ast\ast))(x_k-y_k) + \sum_{j=k-m}^{k-1} \alpha^k_j(I-\eta D^2f(\ast_j))(x_j-x_k),
    \end{split}
\end{equation*}
where we define
\begin{equation}\label{ddf}
\begin{aligned}
        &D^2f(\ast_j)=\int_0^1 D^2f(tx_{k-j}+(1-t)x_k)dt,\\
        &D^2f(\ast\ast)=\int_0^1 D^2f(tx_{k}+(1-t)y_k)dt.
\end{aligned}    
\end{equation}

Note that 
$$
x_k-y_k = x_k - \alpha^k_kx_k - \sum_{j=k-m}^{k-1} \alpha^k_jx_j = \sum_{j=k-m}^{k-1} \alpha^k_j(x_k-x_j).
$$
Hence $b_k$ can be further written as 
\begin{equation}\label{bk}
\begin{split}
b_k &= \sum_{j=k-m}^{k-1} \alpha^k_j(I-\eta D^2f(\ast\ast))(x_k-x_j) + \sum_{j=k-m}^{k-1} \alpha^k_j(\eta D^2f(\ast_j)-I)(x_k-x_j)\\
&= \eta\sum_{j=k-m}^{k-1} \alpha^k_j(D^2f(\ast_j)- D^2f(\ast\ast))(x_k-x_j). 
\end{split}
\end{equation}
Note that $b_k$ would vanish if $f$ were a quadratic function. By mean value theorem, there exists $\omega\in C[0,\infty)$ and is increasing with $\omega(0)=0$ such that
\begin{align*}
\|D^2f(*_j)-D^2f(**)\|
&\leq\omega(\|*_j-**\|)\to 0\quad \text{as}\quad k\to\infty,
\end{align*}
since $x_k,y_k\to x^*$ as $k\to\infty$ implies $*_j=t_1x_{k-j}+(1-t_1)x_k\to x^*$ and $**=t_2x_{k}+(1-t_2)y_k\to x^*$ for some $t_1, t_2\in[0,1]$. Denote $o(k) = \max_{k-m\leq j\leq k-1} \|D^2f(*_j)-D^2f(**)\|$, and assume $x_{j}\in B(x^*, \epsilon)$ with $\epsilon <1$ for $j=k-m, \cdots k$, (\ref{bk}) leads to 
\begin{equation}
\begin{split}
\|b_k\|&\leq o(k)\eta \left\|\sum_{j=k-m}^{k-1} \alpha^k_j(x_k-x_j)\right\|\\
&\leq o(k)\eta \left(\sum_{j=k-m}^{k-1} (\alpha^k_j)^2\right)^{1/2}\left(\sum_{j=k-m}^{k-1}\|x_k-x_j\|^2\right)^{1/2}\\
&\leq o(k)\eta\|\alpha^k_*\|\sqrt{2m}\\
&\leq o(k)\eta\sqrt{2m}\|(U_k^\top U_k+\lambda I)^{-1}U_k^\top\|\|g_k\|,
\end{split}
\end{equation}
where (\ref{lsalpha}) is used in the last inequality. With this bound and \eqref{bda}, \eqref{gab} becomes
\begin{equation*}
    \begin{split}
        \eta\|g_{k+1}\|&\leq \|a_k\|+\|b_k\|, \\
    &\leq\eta(1-\eta\mu)\delta_k\|g_k\|+o(k)\eta c\|g_k\|.
    \end{split}
\end{equation*}
with $ c \geq \sqrt{2m}\|(U_k^\top U_k+\lambda I)^{-1}U_k^\top\|$. Hence
$$
\|g_{k+1}\|\leq (1-\eta\mu)\delta_k\|g_k\|+o(k)c\|g_k\| \\
        =\left((1-\eta\mu)\delta_k+o(k)c\right)\|g_k\|,
$$
which implies
\begin{equation}
       \frac{\|g_{k+1}\|}{\|g_k\|}\leq (1-\eta\mu)\delta_k +co(k)
        =(1-\eta\mu)\delta_k+o(k),
\end{equation}
where
$$
\delta_k=\frac{\|\Pi_k g_k\|}{\|g_k\|}\leq 1.
$$

\section*{Acknowledgments}
This research was partially supported by the National Science Foundation under Grant DMS1812666.

\medskip

\bibliographystyle{amsplain}
\bibliography{ref}

\end{document}